\documentclass[graybox]{svmult}

\usepackage{type1cm}        
\usepackage{makeidx}         
\usepackage{graphicx}        
\usepackage{multicol}       
\usepackage[bottom]{footmisc}

\usepackage{newtxtext} 
\usepackage{newtxmath} 

\makeindex             

\usepackage{amsmath,amssymb,latexsym,color}
\usepackage{enumitem}
\usepackage{hyperref}
\usepackage{setspace}
\usepackage{hyperref}
\usepackage{color}
\usepackage{comment}
\usepackage[utf8]{inputenc}
\usepackage{algorithm,algorithmic}
\usepackage{graphicx}
\usepackage{todonotes}
\usepackage{xargs}
\usepackage{float}
\restylefloat{table}
\newcommandx{\unsure}[2][1=]{\todo[linecolor=magenta,backgroundcolor=magenta!25,bordercolor=magenta,#1]{#2}}
\newcommandx{\suggestion}[2][1=]{\todo[linecolor=blue,backgroundcolor=blue!25,bordercolor=blue,#1]{#2}}

\usepackage[margin = 4 cm]{geometry}

\def\Z{{\mathbb Z}}

\def\G{{\mathcal G}}

\def\K{{\mathcal K}}
\def\S{{\mathcal S}}

\newcommand {\mm}[1] {\ifmmode{#1}\else{\mbox{\(#1\)}}\fi}

\newcommand{\Xcal}        {\mm{\mathcal X}}
\newcommand{\Ycal}        {\mm{\mathcal Y}}

\newcommand{\para}[1]        {\vspace{1mm}\noindent{\textbf{#1}}}

\newcommand{\etal} {\textit{et. al.}}

\newcommand{\blue}[1]{#1}

\begin{document}

\title*{Graph Pseudometrics from a Topological Point of View}
\titlerunning{Graph Pseudometrics from a Topological Point of View} 

\author{Ana Lucia Garcia-Pulido, Kathryn Hess, Jane Tan, Katharine Turner, Bei Wang, Naya Yerolemou}
\authorrunning{AL Garc\'ia-Pulido, K Hess, J Tan, K Turner, B Wang, N Yerolemou}

\institute{A. L. Garc\'ia-Pulido \at Department Of Computer Science, University of Liverpool, Liverpool L69 3BX, UK,\\ \email{A.L.Garcia-Pulido@liverpool.ac.uk}
\and K. Hess \at Laboratory for Topology and Neuroscience, Brain Mind Institute, Ecole polytechnique f\'ed\'erale de Lausanne (EPFL), 1015 Lausanne, Switzerland, \email{kathryn.hess@epfl.ch}
\and J. Tan, N. Yerolemou \at Mathematical Institute, University of Oxford, Oxford OX2 6GG, UK,\, \email{jane.tan@maths.ox.ac.uk, Naya.Yerolemou@maths.ox.ac.uk}
\and K. Turner \at Mathematical Sciences Institute, Australian National University, Canberra, ACT 2601, Australia, \email{katharine.turner@anu.edu.au}
\and B. Wang \at School of Computing, University of Utah, Salt Lake City, UT 84112, USA, \email{beiwang@sci.utah.edu}
}

\maketitle

\abstract{
We explore pseudometrics for directed graphs in order to better understand their topological properties.
The directed flag complex associated to a directed graph provides a useful bridge between network science and topology.  Indeed, it has often been observed that phenomena exhibited by real-world networks reflect the topology of their flag complexes, as measured, for example, by Betti numbers or simplex counts. 
As it is often computationally expensive (or even unfeasible) to determine such topological features exactly, it would be extremely valuable to have pseudometrics on the set of directed graphs that can both detect the topological differences and be computed efficiently. 

To facilitate work in this direction, we introduce methods to measure how well a graph pseudometric captures the topology of a directed graph. We then use these methods to evaluate some well-established pseudometrics, using test data drawn from several families of random graphs.

}

\newpage
\section{Introduction}
\label{sec:introduction}

A typical strategy for studying complex networks is to extract features (i.e., parameters, properties, etc.) of the networks that are simpler to analyze and compare than the networks themselves, yet still capture their essential structures. 
Depending on the context, these features can be local or global and include, for instance, the number of connected components, the number of cycles, the lengths of shortest paths/cycles, \blue{graphlet counts}, and the graph spectrum, as well as many non-numerical properties.

One method for extracting properties of a directed network begins with the construction of its associated \emph{directed flag complex}, which is a topological space built from the directed cliques in the graph. 
Topological features of the directed flag complex provide revealing properties of the original network. 
For instance, the \emph{homology groups} of the directed flag complex reflect how cliques assemble to form the graph globally, enabling us to distinguish between some graphs: if the homology groups of two graphs are different, then the two graphs are not isomorphic. 
From the homology groups one can extract \emph{Betti numbers}, which measure the rank of the homology groups. 
While the $0$-th Betti number is simply the number of connected components of the graph, higher Betti numbers measure how intricately higher dimensional cliques intersect. In contrast to classical graph parameters, these invariants capture higher order structures.

Applications of topological methods to the analysis of networks have been motivated, in particular, by the desire to understand the relation between function and structure of biological networks (see e.g. \cite{Chad2015,Giusti16,Masulli16,ReimannNolteScolamiero2017,Sizemore16,Sizemore18}). 
The work in these articles strongly suggests that biological function reflects topological structure, as measured by Betti numbers for example.
A significant barrier to confirming and exploiting this observation, however, is the computational complexity of determining these features for real-world networks. 
For instance, directly comparing the homology groups or the Betti numbers of two flag complexes is an NP-hard problem~\cite{Adamaszek16}.

In this paper, we explore the hypothesis that there are pseudometrics on the set of directed graphs that detect differences in their topological features. 
As a first step to identifying or constructing such a pseudometric, we introduce methods to measure \blue{how well a given pseudometric reflects differences in the topological features of two graphs.}

We implement our methods of comparison and apply them to existing candidate pseudometrics on the set of directed graphs. 
There are already many high-performing pseudometrics on the space of graphs that take into account 1-dimensional structural features (see, e.g., \cite{Kriege20,Nikolentzos19,Tantardini19}), but there is a dearth of literature on whether they also capture any high-dimensional structure. We note that these methods can also be applied to undirected graphs by considering instead the usual flag complex. However, the analysis in the present paper is restricted to the directed case, since many real-world networks (particularly biological networks) are naturally directed. Our test data set is drawn from several families of directed random graphs for which we can control the parameters and behavior.

We present here an experimental study, comparing topological pseudometrics based on Betti numbers and simplex counts with several well-established pseudometrics for directed graphs. Our analysis is four-fold. First, we study the similarities between clusterings based on these pseudometrics by computing both their Fowlkes-Mallows indices and the distance correlations between them. 
Almost all of the pseudometrics tested are shown to be closely related with high Fowlkes-Mallows indices and high distance correlations. 
Second, we apply \blue{$k$-nearest neighbors} ($k$-NN) classification as a measure of classification accuracy for three models of random graphs. We find that all pseudometrics achieve near perfect classification accuracy. 
Third, we test for relationships between each pseudometric and the various random graph parameters by performing permutation tests with distance correlation and the Fowlkes-Mallows indices. 
Using the permutation tests, we can reject the null hypothesis of independence between all the different pseudometrics when considered over pooled sets of directed graphs with multiple parameter values. However, we cannot in general reject the null hypothesis of independence once we restrict to a specific model and parameter value. This indicates that the latent variable of the parameter of the model is important.  
Finally, we apply $k$-NN regression to our pseudometrics to try to predict the topological feature vectors of given graphs. 


\para{Outline} We begin by reviewing the requisite topological and combinatorial background in Section~\ref{sec:prelim}. This includes the definitions of our topological feature vectors together with the topological pseudometrics that they induce, as well as a brief introduction to each of the existing pseudometrics that we compare to our topology-based ones: TriadEuclid, TriadEMD, and Portrait Divergence. 
We present two methods for comparing pseudometrics in Section~\ref{sec:compare}, one based on clustering and the other on distance correlation. Here, we also describe the permutation test and $k$-NN regression, as well as the classification methods that we use. Technical details of our experiments can be found in Section~\ref{sec:methods}. The final results of the comparison are presented in Section~\ref{sec:results}.

\section{Directed Graph Pseudometrics}
\label{sec:prelim}

A \emph{directed graph} $G$ consists of a set of \emph{vertices} $V$ together with a set of \emph{edges} $E$, which are ordered pairs of vertices. All graphs in this paper are \emph{finite}, meaning $V$ and $E$ are both finite sets.
The direction of an edge $(u,v)\in E$ is taken to be from $u$, the {\em origin}, to $v$, the {\em destination}.
We require that $G$ not contain any self-loops, that is, for each $(u,v)\in E$, $u \neq v$. Secondly, for each pair of vertices $(u,v) \in V \times V$, there is at most one directed edge from $u$ to $v$.
Note, however, that we do allow both $(u,v) \in E$ and $(v, u) \in E$. In other words, the directed graphs we consider are simple except for bigons, and we shall simply refer to them as \emph{graphs} or \emph{digraphs}.

We use several standard definitions associated with digraphs. 
The \emph{out-degree} of a vertex $v$ is the number of edges having $v$ as the origin, while its \emph{in-degree} is the number of edges having $v$ as the destination. 
A vertex $v$ is a \emph{sink} if its out-degree is zero and its in-degree is at least one and  a \emph{source} if its in-degree is zero and its out-degree is at least one. 
A \emph{path} in $G$ is a list of distinct edges such that the destination of the $i$-th edge is the same as the origin of the $(i+1)$-st edge and such that no vertex is traversed more than once, except the path may end at the vertex where it started.  If the first and last vertices of a path are the same, we call it a \emph{cycle}. 

Let $\G$ denote the set of all finite directed graphs. We will view $\G$ as a space endowed with several natural \emph{pseudometrics}. First, recall that a pseudometric on a set $\Xcal$ is a function $d_\Xcal \colon \Xcal \times \Xcal  \to [0,\infty)$ such that for all $x,y,z\in X$,
\begin{enumerate}
        \item $d_\Xcal(x,x) = 0$;
        \item $d_\Xcal(x,y) = d_\Xcal(y,x)$ (symmetry); and
        \item $d_\Xcal(x,z) \leq d_\Xcal(x,y) + d_\Xcal(y,z)$ (triangle inequality).
\end{enumerate}
Importantly, points need not all be distinguishable by a pseudometric: it is possible that $x \neq y$, even though $d_\Xcal(x,y)=0$. The pair $(\Xcal, d_\Xcal)$ is a \emph{pseudometric space}. 

We will describe three well-established pseudometrics on $\G$: TriadEuclid in Section~\ref{sec:triadEuclid}, TriadEMD in Section~\ref{sec:triadEMD}, and portrait divergence in Section~\ref{sec:portrait-divergence}. Before this though, we define two topological summaries for elements in $\G$ based on Betti numbers and simplex counts in Section~\ref{sec:betti-simplex}, which provide a crucial point of comparison.


\subsection{Betti Numbers and Simplex Counts}
\label{sec:betti-simplex}

The key construction we study in this paper is the directed flag complex of a directed graph. For a more detailed account we refer to the work of Luetgehetmann \etal~\cite{LuetgehetmannGovcSmith2019} and Reimann \etal~\cite{ReimannNolteScolamiero2017}. We assume familiarity with standard notions of abstract simplicial complexes and simplicial homology, introduced for instance in~\cite{Hatcher2002,Munkres1984}.

\begin{definition}[Abstract Directed Simplicial Complex]
An \emph{abstract directed simplicial complex} on a vertex set $V$ is a collection $\K$ of lists (i.e., totally ordered sets) of elements of $V$ such that for every sequence $\sigma \in \K$, every subsequence $\tau$ of $\sigma$ belongs to $\K$. 
\end{definition}

An element $\sigma \in \K$ is called a \emph{(directed) simplex}.
If $\sigma$ is of length $p+1$, then we call it an \emph{$p$-simplex}.  
The collection of $p$-simplices of $\K$ is denoted  $\K_p$.
If $\sigma\in \K$ and $\tau \subset \sigma$, then $\tau$ is  called a \emph{face} of $\sigma$.
The  $i$-th face of an $p$-simplex $\sigma = (v_0,\ldots,v_p) $ is the
$(p-1)$-simplex obtained by removing the  $v_{i}$ from the list $\sigma$.

The notion of an abstract directed simplicial complex is a variant of the more common notion of abstract simplicial complex. Henceforth, we always mean abstract directed simplicial complexes when we say \emph{simplicial complexes}. The following definition illustrates how directed simplicial complexes arise naturally from directed graphs.

\begin{definition}[Directed Flag Complex]
Given a directed graph $G = (V,E)$, the {\em directed flag complex} of $G$, denoted $\K = \K(G)$, is defined as follows. 
    \begin{itemize}    
        \item Take $\K_0 = V$.
        \item For $p \geq 1$, a {\em directed p-simplex} $\sigma$ in $\K_n$ is an $(p+1)$-tuple of vertices $(v_0, \ldots, v_p)$ such that there is a directed edge from $v_i$ to $v_j$ for every pair of vertices $v_i$, $v_j$ with  $0 \leq i<j\leq p$.
    \end{itemize}
\end{definition}    
For $\sigma = (v_0, \ldots, v_p) \in K$, we call $v_0$ the {\em source} of the simplex, since there is a directed edge from $v_0$ to $v_i$ for every $i>0$. 
Similarly, we call $v_p$ the {\em sink} of the simplex, since there is a directed edge from $v_i$ to $v_p$ for every $i < p$. Note that these are consistent with the equivalent digraph notions \blue{since a $p$-simplex is characterised by the \emph{ordered} sequence of vertices and not by the underlying set of vertices.}

Throughout this paper,  we consider simplicial homology of directed flag complexes with $\mathbb{F}_2$-coefficients, where homology is defined in the usual way. 
For a simplicial complex $\K$, let $H_p(\K)$ denote its \emph{$p$-th homology group} and $\beta_p(\K)=\dim H_p(\K)$ its \emph{$p$-th Betti number}. 
        
Given $G \in \G$, we work with two simple topology-based feature vectors on $G$, defined as follows. 
Let $p \in \Z_{\geq 0}$ denote the maximum homological dimension of interest, which is context-dependent; in our experiments, $p = 6$. 
For $0\leq k\leq p$, let $b_k(G) = \max \{0, \log\left(\beta_k\left(\K\left(G\right)\right)\right)\}$ and 
$$b(G) = \big(b_0(G), \dots, b_p(G)\big).$$
Let $\gamma_k(G)$ be the number of directed $k$-simplices of $\K(G)$, $c_k(G)= \max \{0, \log\left(\gamma_k\left(\K\left(G\right)\right)\right)\}$, and
$$c(G) = \big(c_0(G), \dots, c_p(G)\big).$$  
The vectors $b(G)$ and $c(G)$ consist of the logarithms of  the Betti numbers and simplex counts of $\K(G)$, respectively. 
Let $\| \cdot \|_2$ denote the Euclidean norm on $\mathbb R^n$.

\smallskip
\begin{definition}[Topological Pseudometrics] The pseudometrics $d_\beta$ and $d_\Delta$ on $\G$ are specified by 
\begin{align}
& d_{\beta}(G,G') = \| b(G) - b(G')\|_{2};\\
& d_{\Delta}(G,G') = \| c(G) - c(G')\|_{2}.   
\end{align}
for any pair of directed graphs $G, G' \in \G$.
\end{definition}

\subsection{TriadEuclid} 
\label{sec:triadEuclid}
   
\blue{First introduced by Przulj \etal~\cite{PrzuljCorneilJurisica2004}}, the term {\em graphlet} is often used to mean a small connected graph up to a fixed size. A {\em graphlet-based} pseudometric compares the graphlet counts in pairs of graphs. Graphlet-based pseudometrics are a well established tool in network analysis, as they capture local structural similarity of two graphs. 

Xu and Reinert~\cite{XuReinert2018} introduced two graphlet-based pseudometrics that outperform the best previously defined directed graphlet methods in graph classification tasks: {\em TriadEuclid} and {\em TriadEMD}. Both of these consider only directed graphlets on three vertices. In this section we focus on TriadEuclid, which  measures the difference between $3$-graphlet counts of two directed graphs in terms of their Euclidean distance. 
We remark that there are exactly $13$ isomorphism classes of connected directed graphlets with $3$ vertices (referred to as $3$-graphlets). A complete list can be found in \cite{XuReinert2018}.

Let $G=(V,E)$ a directed graph. If the induced subgraph determined by three vertices of $G$ is connected, it must be isomorphic to one of the $13$   $3$-graphlets. Let $n_i(G)$ be the number of induced subgraphs of $G$ that are isomorphic to the $i$-th graphlet, $i\in I=\{1,\ldots,13\}$.  Define $\phi(G)\in \mathbb N^{13}$ by
\begin{align}
\phi(G)_j=\frac{n_j(G)}{\sum_{i\in I}n_i(G)}
\end{align}
for $j\in I$.
\begin{definition} The {\em TriadEuclid} pseudometric on the set $\G$ of finite directed graphs
is defined by
\begin{align}
\operatorname{TriadEuclid}(G,G') = \|\phi(G)-\phi(G')\|_2,
\end{align} 
for all $G,H\in \G$, where $\|\cdot\|_2$ denotes the Euclidean norm.
\end{definition}

The complexity of this algorithm is $O(nd^2)$, where $d$ is the maximum degree of vertices in $G$ and $G'$, $n$ is the maximum number of vertices in $G$ and $G'$.

\subsection{TriadEMD}
\label{sec:triadEMD}

TriadEMD is another graphlet-based pseudometric defined by Xu and Reinert \cite{XuReinert2018}, computed in terms of the earth mover distance between generalized degree distributions of two directed graphs.

\begin{definition}
Given a directed graph $H= (V, E)$, an {\em automorphism} of $H$ is a graph isomorphism $h\colon H\to H$. The set of automorphisms of $H$ forms a group under composition,  denoted $\operatorname{Aut}(H)$, which acts on $H$.

Let $v\in V$. The {\em orbit} of $v$ under the action of $\operatorname{Aut}(H)$ is defined by
\begin{align}
\operatorname{Orb}(v) = \{h(v)| h\in \operatorname{Aut}(H)\}.
\end{align}
\end{definition} 

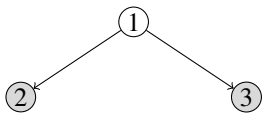
\begin{figure}
    \centering
    \begin{tikzpicture}
    \node[circle,draw=black, fill=white, inner sep=1pt,minimum size=5pt] (1) at (1,1) {1}; 
    \node[circle,draw=black, fill=gray!30!white, inner sep=1pt,minimum size=5pt] (2) at (-0.5,0) {2}; 
    \node[circle,draw=black, fill=gray!30!white, inner sep=1pt,minimum size=5pt] (3) at (2.5,0) {3}; 
    \draw [->] (1) -- (3);
    \draw [<-] (2) -- (1);
    \end{tikzpicture}
    \caption{A graph with two vertex orbits $\{1\}$ and $\{2,3\}$.}
    \label{fig:vertex_orbit}
\end{figure}

Projecting to orbits, one considers the genuinely different positions that a vertex takes in a fixed graphlet. 
\blue{~\autoref{fig:vertex_orbit} shows an example of a three-vertex graph with only two distinct orbits:
\[
\operatorname{Orb}(1) = \{1\},\quad
\operatorname{Orb}(2) = \operatorname{Orb}(3) =\{2,3\}.
\]}

A complete list of the $30$ orbits of $3$-graphlets can be found in \cite{XuReinert2018}. A list of graphlets with between two and four vertices is given in \cite{sarajlic16}.

Degree distributions can be used to study {\em hubs} in a network, i.e., vertices of high degree. 
The {\em degree distribution} $P$ of an undirected graph is defined so that $P(k)$ is the proportion of vertices in $G$ with degree $k$:
\begin{align}
P(k)= \frac{\#\{v\in V| \operatorname{deg}(v) = k\}}{\# V}.
\end{align}
When $G$ is directed, there are analogous in-degree and out-degree distributions. 

TriadEMD is defined in terms of the degree distributions of orbits in triads. 
For any orbit $i$ \blue{of a fixed graphlet}, the \emph{orbit-i-degree} of a vertex $v$ of a digraph $G$  is the number of copies  of orbit $i$ in $G$ of which $v$ is a vertex~\cite{XuReinert2018}. 
The {\em orbit-$i$-degree distribution} $P_i$ of $G$ is defined so that $P_i(k)$ is the proportion of vertices in $G$ with orbit-$i$-degree $k$. 

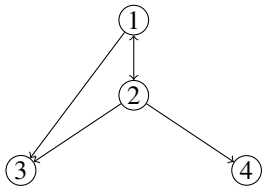
\begin{figure}
    \centering
    \begin{tikzpicture}
    \node[circle,draw=black, fill=white, inner sep=1pt,minimum size=5pt] (1) at (1,2) {1};
    \node[circle,draw=black, fill=white, inner sep=1pt,minimum size=5pt] (2) at (1,1) {2}; 
    \node[circle,draw=black, fill=white, inner sep=1pt,minimum size=5pt] (3) at (-0.5,0) {3}; 
    \node[circle,draw=black, fill=white, inner sep=1pt,minimum size=5pt] (4) at (2.5,0) {4}; 
    \draw[<->] (1) -- (2);
    \draw[->] (1) -- (3);
    \draw[->] (2) -- (3);
    \draw[->] (2) -- (4);
    \end{tikzpicture}
    \caption{Graph $G$ with $\operatorname{Orb}\big(\{1\}\big)$-degree distribution given in \autoref{table:orbit_distribution}.}
    \label{fig:orbit_distribution}
\end{figure}

\blue{Consider, for example, the graph $G$ on four vertices shown in ~\autoref{fig:orbit_distribution}, for which we now calculate the $\operatorname{Orb}\big(\{1\}\big)$-degree distribution of the graphlet in ~\autoref{fig:vertex_orbit}. The $\operatorname{Orb}\big(\{1\}\big)$-degrees of vertices $1$, $2$, $3$, and $4$ are $1$, $3$, $0$, and $0$ respectively; ~\autoref{table:orbit_distribution} gives the $\operatorname{Orb}\big(\{1\}\big)$-degree distribution.}

\begin{table}[H]
\centering
\begin{tabular}{ |c | c | c | c | c| }
 $k =$degree & 0 & 1 & 2 & 3 \\ 
 \hline
 $P_{\{1\}}(k)$ & 0.5  & 0.25 & 0 & 0.25\\  
\end{tabular}
\caption{$\operatorname{Orb}\big(\{1\}\big)$-degree distribution of $G$ ~\autoref{fig:orbit_distribution}}
\label{table:orbit_distribution}
\end{table}

The original definitions of the orbit degrees and orbit graphlet metrics for undirected and directed graphs, together with additional examples can be found in \cite{Przulj10,sarajlic16,wegner18,XuReinert2018}.

\begin{definition} 
Let $\Theta$  be the set of 30 orbits of $3$-graphlets. The \emph{TriadEMD} pseudometric on the set $\G$ of directed graphs is given by
\begin{align}
    \operatorname{TriadEMD}(G,G') = \frac{1}{30}\sum_{i\in \Theta} \operatorname{EMD}(P_i,P'_i),
\end{align} 
where $P_i$ (resp. $P'_i$) denotes the orbit-$i$-degree distribution of $G$ (resp. $G'$), and $\operatorname{EMD}$ denotes the earth mover distance (EMD) between the distributions. Recall that the EMD between two probability distributions $P$ and $Q$ on the real line is defined by
\begin{align}
\operatorname{EMD}(P,Q)=\int_{-\infty}^{\infty}|F_P(x) - F_Q(x)|dx,
\end{align}
where $F_P$ and $F_Q$ are the cumulative distribution functions of $P$ and $Q$ respectively. 
\end{definition}
The complexity of this algorithm is $O(n \log n)$, where $n$ is the maximum between the numbers of vertices of $G$ and of $G'$.

\subsection{Portrait Divergence}
\label{sec:portrait-divergence}

First introduced by Bagrow and Bollt \cite{BagrowBollt2019}, {\em Network Portrait Divergence} is a pseudometric  that compares the distributions of the shortest-path lengths of two graphs.  \blue{Recall that the {\em diameter} of a finite directed graph is the maximum over all pairs of vertices of the length of the shortest (directed) path from the first vertex to the second.}

\begin{definition}
Let $G=(V,E)$ be a directed graph of diameter $d$, with $n$ vertices, and denote by $d_{sp}(v,u)$ the shortest path distance from $v$ to $u$. For every vertex $v\in V$ and $0\leq l\leq d$, define
\begin{align}
s^l_v =\# \{ u\in V| d_{sp}(v,u) = l\}.    
\end{align}
The {\em network portrait} \cite{Bagrow08}
of $G$ is the matrix $B = \left(B_{lk} \right)$ with entries 
\begin{align}
    B_{lk} = \#\{v\in V| s^l_v = k\},
\end{align}
where $0\leq l\leq d$ and $1\leq k \leq n-1$. That is,  $B_{lk}$ is the number of vertices that are at distance $l$ from exactly $k$ other vertices.
\end{definition}
The first row of the network portrait $B$ of $G$ is precisely the degree distribution of $G$.
The second row is the degree distribution of next-nearest neighbors, and so on. 
The network portrait $B$ also captures structural features of $G$ such as the number of edges, the diameter of $G$, and the distribution of shortest paths. 
It has been shown to be a graph invariant~\cite{BagrowBollt2019}. 
The complexity of computing an unweighted portrait is $O(mN+N^2)$, due to the procedure of finding minimum-length paths, where $m$ is the number of edges and $N$ the number of vertices in $G$. 

Let $G= (V,E)$ be a directed graph. 
Let $P(k,l)$ be the probability of \blue{choosing two vertices $(u, v)\in V\times V$ uniformly at random such that $d_{sp}(u,v) = l$ and $s^l_u = k$. If $s^l_u = k$, then there exist $k$ vertices $u_1,\ldots,u_k$ with $d_{sp}(u,u_i) = l$. Therefore,
\begin{align*}
    \# \{ (u, v)\in V\times V: d_{sp}(u,v) = l , s^l_u = k \} 
    &= k\cdot\#\{v\in V: s^l_v = k\} \\
    &= k\cdot B_{lk},
\end{align*} and in consequence,
\begin{align}
P(k,l) = \frac{k B_{lk}}{n^2}.
\end{align}
}
Given graphs $G$ and $G'$, we define (joint) distributions $P$ and $P'$ for all rows of their portraits $B$ and $B'$.
The KL-divergence between $P$ and $P'$ is then defined by 
\begin{align}
KL(P || P') = \sum_{l=0}^{\max(d,d')} \sum_{k=0}^N P(k,l) \log \frac{P(k,l)}{P'(k,l)}. 
\end{align}


\begin{definition}
The {\em Portrait Divergence} (PD) pseudometric on the set $\G$ of finite directed graphs is the Jensen-Shannon divergence,  
\begin{equation}
PD(G,G') = \frac{1}{2}KL(P||M) + \frac{1}{2}KL(P'||M),
\end{equation}  for all $G,G'\in \G$, where $M=(P+P')/2$ is the mixture distribution of $P$ and $P'$. 
\end{definition}

\section{Statistical Tools}
\label{sec:compare}

We are interested in comparing known pseudometrics on $\G$ (e.g., TriadEuclid, TriadEMD, and PD) to $d_\beta$ and $d_\Delta$, which requires some care. A natural first idea would be to use Gromov-Hausdorff distance, but since it  measures how close two spaces are to being isometric, it is too rigid for practical purposes. For example, changing even just one Betti number of a single graph can result in a very large Gromov-Hausdorff distance. 
Multiplying a topological pseudometric by a constant has a similar effect as well. 

Instead, we employ two comparison methods that are invariant under rescaling either pseudometric and more robust under perturbations of points. The first method is based on the distance correlation between pseudometrics (Section~\ref{sec:dist_cor}), while the second is based on the Fowlkes-Mallows index of the clusterings of the two pseudometrics (Section~\ref{sec:clustering_comparison}). Both rely on taking finite samples $\S \subset \G$ and performing an analysis on $\S$.

We remark that similar comparison methods can also be applied to undirected graphs by replacing the directed flag complex with the flag complex, though we restrict to the directed case in this paper.

\subsection{Distance Correlation}
\label{sec:dist_cor}

\blue{The concept of distance correlation was first introduced in by Szekely \etal\cite{SzekelyRizzoBakirov2007} for two paired random vectors in Euclidean space and generalized to metric spaces by Lyons~\cite{Lyons2013}. Distance correlation measures linear and non linear relationships between two distributions lying in possibly different metric spaces. We follow the definition of the sample distance correlation of a paired sample as formulated by Turner and Spreemann~\cite{Turner19}. We use the sample distance correlation as an estimation of the distance correlation and refer to the sample distance correlation simply as the distance correlation.}

Let $(\Xcal,d_{\Xcal})$ and $(\Ycal,d_{\Ycal})$ be pseudometric spaces, and let $(X,Y) = \big\{(x_i,y_i)\big\}_{1\leq i \leq l} \subset \Xcal\times \Ycal$ be paired samples. 
For $1\leq i,j\leq l$, let $a_{i,j} = d_{\Xcal}(x_i,x_j)$ and $b_{i,j} = d_{\Ycal}(y_i,y_j)$, so that $a = (a_{i,j})$ and $b = (b_{i,j})$ denote  matrices of pairwise distances in $\Xcal$, and $\Ycal$, respectively. Let $\overline{a}^i$ and $\overline{b}^i$ denote the row means and $\overline{a}_j$ and $\overline{b}_j$ the column means of the matrices $a$ and $b$. Let $\overline{a}$ and $\overline{b}$ denote the total matrix means. 
Define {\em doubly centred matrices} $(A_{k,l})$ and $(B_{k,l})$ by $A_{k,l}=a_{k,l}-\overline{a}^k -\overline{a}_l +\overline{a}$ and $B_{k,l}=b_{k,l}-\overline{b}^k -\overline{b}_l +\overline{b}$. 
\begin{definition}
    The {\em sample distance covariance} of the paired sample $(X,Y)$ is
 \begin{align}
    \operatorname{dcov}(X,Y)=\frac{1}{n^2}\sum_{k,l=1}^n A_{k,l} B_{k,l}.
  \end{align}
  \blue{When $\operatorname{dcov}(X,Y) \geq 0$, let $\operatorname{dCov}(X,Y) = (\operatorname{dcov}(X,Y))^{1/2}$.}
  
  The {\em sample variance} of the sample $X$ is defined to be
 \begin{align}
    \operatorname{dVar}(X)=\left(\frac{1}{n^2}\sum_{k,l=1}^n A_{k,l}^2\right)^{1/2}.
\end{align}
  If $\operatorname{dVar}(X)\operatorname{dVar}(Y)\neq 0$ and $\operatorname{dcov}(X,Y) \geq 0$, the {\em sample distance correlation} is given by
  \begin{align}
      \operatorname{dCor}(X,Y)= \frac{\operatorname{dCov}(X,Y)}{(\operatorname{dVar}(X)\operatorname{dVar}(Y))^{1/2}}.
  \end{align}
  If $\operatorname{dVar}(X)\operatorname{dVar}(Y)= 0$, then we set $\operatorname{dCor}(X,Y)=0$.
\end{definition}
\blue{For strongly negative metrics spaces, which includes Euclidean space, the sample distance covariance between $X$ and $Y$ is always non-negative, and $X$ and $Y$ are independent if and only if $\operatorname{dCov}(X,Y) = 0$, see~\cite{Lyons2013}.} This is not true for all metric spaces, but in all of our cases the sample covariance is non-negative, and therefore the sample distance correlation is well-defined.

Notice that sample distance correlation is invariant under scalar multiplication of either of the two metrics. Cauchy-Schwarz implies that $|\operatorname{dCor}(X,Y)|\leq 1$ and that equality is attained when the doubly centred matrices are scalar multiples of each other. Thus, high correlation measures the (possibly non-linear) relationship between $X$ and $Y$ arising from a linear relationship between their corresponding metrics.

\blue{In our case we have $\Ycal = \Xcal$, and we consider the sample distance correlation on  paired samples $(X,X)$. }

\subsection{Fowlkes-Mallows Index}
\label{sec:clustering_comparison}

Our second method of comparing pseudometric spaces is based on clustering. Suppose we have pseudometric spaces $(\Xcal,d_0)$ and $(\Xcal,d)$. Given a finite sample $X \subset \Xcal$, we can compute hierarchical clusterings $A_0$ and $A$ corresponding to both pseudometrics and hence apply standard techniques of cluster evaluation. In particular, we use the Fowlkes-Mallows index~\cite{FowlkesMallows1983}, which provides a measure of similarity between two clusterings. Making such comparisons across many choices of sample $X$ provides a means of comparing $(\Xcal,d)$ to $(\Xcal,d_0)$.

\blue{
We employ a standard agglomerative hierarchical clustering algorithm (see documentation of~\cite{scipycluster} for details), initialised with every point in a separate cluster. In each step, two clusters of minimum distance to each other are chosen and merged into one to move up the hierarchy. We use complete linkage, so the distance between clusters is the maximum distance between any two points of those clusters. The algorithm terminates when all points are in a single cluster.} 

\blue{
In some applications, the number of expected clusters is known from the outset, allowing early termination of the algorithm. As this does not apply to our case, we first run the clustering algorithm and then use \emph{silhouette analysis} to choose the most natural number of clusters with respect to $d_0$. This method, introduced by Rousseeuw~\cite{Rousseeuw1987}, assesses how well a clustering captures the structure of the data using only internal distance information. Informally, we judge how well a single point $x$ has been placed within a cluster based on a {\em silhouette value}, which quantifies how close $x$ is to the other points in its own cluster in contrast to points in other clusters, and then extend to a numerical score for the clustering.}

\begin{definition}
Suppose that we have points in a pseudometric space $(\Xcal,d)$  partitioned into clusters $C_1,\ldots,C_k$. Given a point $u\in C_i$, let $a(u)$ be the mean distance between $u$ and other points in $C_i$, that is 
\begin{align}
a(u) = \frac{1}{|C_i|-1} \sum_{v \neq u, v \in C_i}d(u,v).
\end{align}
For every $j \neq i$, let $a_j(u)$ be the mean distance between $u$ and points in cluster $C_j$,
\begin{align}
a_j(u) = \frac{1}{|C_j|} \sum_{v \in C_j}d(u,v).
\end{align} 
Set $b(u) = \operatorname{min}_{j\neq i}a_j(u)$. 
If $|C_i| > 1$, define the {\em silhouette value} of $u$ by 
\begin{align}
s(u) = \frac{b(u)-a(u)}{\operatorname{max}\{a(u),b(u)\}},
\end{align}
and set $s(u) = 0$ if $|C_i| = 1$.
\end{definition}

It follows immediately from the definition that $s(u)\in [-1,1]$ and that
\begin{itemize}
    \item $s(u) \approx 1$ indicates that $u$ is appropriately clustered; 
    \item $s(u)\approx -1$ means that $u$ is closer to points from a different cluster,
    \item if $s(u)\approx 0$, then  $u$ is almost equidistant between $C_i$ and at least one other cluster $C_j$.
\end{itemize} 

Based on silhouette values of the points in each of the clusters, one can define a score for the entire clustering.
\begin{definition}
The {\em silhouette coefficient} of the clustering $C_1, \ldots, C_k$ of $l$ points $u_1,..,u_l$ is their  average silhouette value:
\begin{align}
SC = \frac{1}{l}\sum_{1\leq i\leq l} s(u_i).
\end{align}
\end{definition}
The definition above enables us to determine if the clustering $C_1,\ldots, C_k$ reflects a natural structure present in our samples or if instead it seems forced. 

\blue{The actual comparative element in our method based on the silhouette coefficient comes from the use of the Fowlkes-Mallows index \cite{FowlkesMallows1983}, which measures either the similarity between two clusterings with $k$ clusters or the similarity between one clustering and a benchmark classification. We use the former, but both proceed by considering whether pairs of points are consistently assigned to clusters or classified.} Specifically, let $X$ be a set with $n$ elements, and let $A$ and $A'$ be two clusterings of $X$ with $k$ clusters each. Say that a pair of points are \emph{together} in a clustering if and only if they are assigned to the same cluster.
Let $TP$ (true positives) be the number of pairs of points that are together in both $A$ and $A'$, $FP$ (false positives) the number of pairs together in $A$ but not in $A'$, $FN$ (false negatives)  the number of pairs together in $A'$ but not in $A$, and $TN$ (true negatives) the number of pairs that are not together in either $A$ or $A'$.
\begin{definition}  
    The {\em Fowlkes-Mallows} index $FM_k$ of the pair $A, A'$ (for clusterings with $k$ clusters) is given by 
\begin{align}
    \label{eqn:FM_index}
     FM_k := \sqrt{\frac{TP}{TP+FP}\cdot\frac{TP}{TP+FN}}.
 \end{align}
\end{definition}

Notice that $FM_k\in[0,1]$ and that the clusterings $A$ and $A'$ are identical when $FM_k=1$ and the most dissimilar when $FM_k=0$. 

In our experiments, given a directed graph, we apply silhouette analysis to determine the number of clusters $k$ that gives the highest silhouette score for $d_{\beta}$ (respectively, $d_\Delta$). 
We then compare a known pseudometric $d$ (triadEuclid, triadEMD, and PD) to $d_{\beta}$ (respectively, $d_\Delta$) in terms of the Fowlkes-Mallows index of their associated clusterings. 

\para{Remark } A particular challenge in comparing clusterings is the absence of benchmark datasets. It is also difficult to determine the optimal number of clusters for a given dataset. 
In order to circumvent this issue, we use silhouette analysis to determine the most appropriate number of clusters. 
In our experiments, we obtain a unique number of clusters that realizes the maximal silhouette coefficient; however, this is not guaranteed in general. In the case where multiple numbers of clusters attain the same maximum silhouette score, one could easily adapt the present method  by computing Fowlkes-Mallows indices for every viable number of clusters and then taking, for instance, the mean score.

\subsection{Permutation Tests for Paired Data}

The permutation test has become a default method for testing independence. It is a provably valid and consistent test for any consistent dependency measure, such as distance correlation. For more details concerning permutation tests with distance correlation we refer to \cite{ShenPriebeVogelstein2020}.

A permutation test for paired data involves computing some summary statistic (in our case, either distance correlation or the Fowlkes-Mallows index) on the original pairing and then comparing this summary statistic to that of a shuffled pairing.  If the two pseudometrics are truly independent, then the percentile of the summary statistic amongst all those with permuted pairings is drawn uniformly over possible percentiles. 
 
In particular, we can test for a null-hypothesis of independence for two pseudometrics with a permutation test where the sampled \emph{$p$-value} is the percentile of the original distance correlation (resp.\ Fowlkes-Mallows index) among those for the permuted pairings. 

\blue{
Since we compute many $p$-values for comparing various pairs of pseudometrics over many different sample sets of directed graphs, we must correct for multiple hypothesis testing. The \emph{family-wise error rate} is the probability of making one or more false discoveries when performing multiple hypotheses tests. The simplest, most conservative method to bound the  family-wise error rate is the Bonferroni correction, which can be applied to any collection of hypothesis tests, regardless of any dependency structure among the variables. In the Bonferroni correction for $N$ tests, we simply divide the target significance level by $N$; to bound the family-wise error rate by $\alpha$ we use $\alpha/N$ as the threshold for rejecting the null hypothesis. In order to apply the Bonferroni correction, we consider in our analyses each table of $p$-values as a family. 
}

\blue{
To consider all the permutation tests in all of the tables of this paper simultaneously, we must instead bound the \emph{false discovery rate}: the number of discoveries (tests where we reject the null hypothesis) that may be false positives as a proportion of all discoveries made, including both true and false discoveries. To affirm that the expected number of discoveries that are false is at most $\alpha$, we apply the conservative approach by Benjamini and Yekutieli \cite{BenjaminiYekutieli2001} formulated below, since we do not know what dependency structures there may be among the variables,
}

\blue{
\begin{theorem}\cite{BenjaminiYekutieli2001}
Let $C(N)=\sum_{j=1}^N \frac{1}{j}$. Let $p_{(1)}\leq p_{(2)} \leq \ldots \leq p_{(N)}$ be the
ordered observed $p$-values. Define
\begin{align} 
k:=\max\left\{i: p_{(i)}\leq \frac{i}{m C(N)}\alpha \right\}.
\end{align}
If no such $i$ exists, reject no hypothesis.
If an $i$ exists, and we reject the null hypotheses $H_{(1)}, H_{(2)}\ldots H_{(k)}$, then we have controlled the false discovery rate at a level less
than or equal to $\alpha$.
\end{theorem}
}

\subsection{$k$-NN Classification and Regression}

Given a function $f$ on a data set and a set of training data for which the function value is known, one intuitive approach to estimating $f$ on a test point is to extrapolate from the function values on the training points closest to it. 
This is the method of \emph{$k$-nearest neighbors} ($k$-NN). 
There are two slightly different algorithms for this method, depending on whether the unknown function is categorical or real-/vector-valued. 

For a categorical functions, one applies \emph{$k$-NN classification}. 
For a given test point, one considers the $k$ nearest training points and assigns to the test point the category that appears the most frequently among them. We measure success in terms of \emph{classification rate}, i.e., the proportion of classifications that are correct. 

For real- or vector-valued functions,  \emph{$k$-NN regression} is more appropriate. 
For a given test point, one considers the $k$ nearest training points and assigns to the test point the average of their function values. 
The difference between the function estimate and the true function value is called the \emph{residual}. The goal is for the residuals to be  small as possible. We therefore measure success via the \emph{mean squared error} (MSE), i.e., the mean of the squares of the errors/residuals. 

There are many ways to split data into training and testing sets. 
In this paper, we use leave-one-out, a form of cross-validation,  where at each stage we pick a single data point as the testing set and use everything else as the training set. 
We then repeat the process so that every data point is considered as the test sample at some point. 

\section{Experimental Setup}
\label{sec:methods}

In this section, we provide the technical and implementation details of our experimental setup. 
First, we describe the random graph models used to generate our test graphs and the associated parameter selection. 
We then compute topological pseudometrics $d_\beta$ and $d_\Delta$ for the collections of test graphs \blue{in order to perform a comparative study}. 

\subsection{Random Graph Models and Parameters}
\label{sec:random-graphs}

In our experiments, we analyze directed graphs from families based on three random graph models: directed Erd\"{o}s--R\'{e}nyi random graphs (ER), directed geometric random graphs (GR), and random $k$-out graphs with preferential attachment (PA). There are numerous standard references for these models (see, for instance, \cite{Newman18}). \blue{We focused on these three random models to initialise the study of pseudometrics from a topological perspective. In particular, GR and PA are used to model real world networks, such as mobile ad hoc networks, the World Wide Web, and social networks.} We include a brief description with the purpose of describing the specifications used to generate our test data.

\para{Random graph models }
A directed ER random graph is generated by starting with a fixed set of $n$ vertices and adding a directed edge between each ordered pair of vertices independently with probability $\rho$. Note that the edges $u \to v$ and $v \to u$ are also chosen independently of each other, and in particular it is possible for both to be present.

\blue{
A classic geometric random graph is generated by placing vertices uniformly at random in the unit square, and then adding an edge between two vertices whenever the (Euclidean) distance between the vertices is at most equal to a fixed parameter $r$. We consider (biased) oriented GR random graphs obtained by taking an undirected graph generated in the classical sense with vertex set $\{1,\ldots,n\}$, and then for each edge $uv$ (with $u<v$) choosing the direction $u\to v$ with probability 1/3 and $v\to u$ with probability 2/3. The directions are chosen independently for each edge in the undirected graph, but in this collection it is not possible to have both directed edges between a single pair of vertices.}

A PA random graph with parameter $k$ is generated as follows: give each vertex an initial weight of $1$, and select a vertex $u$ with out-degree less than $k$, uniformly at random. 
Choose another vertex $v$ with probability proportional to its weight, add a directed edge from $u$ to $v$ and increase the weight of $v$ by $1$. 
This process terminates when every vertex has out-degree $k$. 
The initial output may have parallel (repeated) directed edges, which we subsequently reduce to a single directed edge leaving at most one edge in each direction between any pair of vertices.

\blue{
\para{Test graphs }
Our test graphs consist of two collections of directed graphs, all on 500 vertices, generated according to the preceding descriptions. The first collection consists of $120$ graphs, with 10 for each of the following parameters:  ER with $\rho\in\{ 0.03,0.06,0.1,0.15,0.2,0.25\}$, GR with $r \in \{0.1,0.175,0.3\}$, and PA with $k \in \{20,40,70\}$. 
We refer to this collection as the {\em point-drawn collection}, since the parameters are chosen from a discrete set of values. }

\blue{
For the second collection, we generated $300$ directed graphs with 100 for each of the three random graph models. The parameters for a fixed model were obtained by generating 25 values independently uniformly at random from a set of four predefined intervals: $\rho$ values are chosen from intervals in $\{(0,0.01), (0.02,0.03), (0.05,0.07), (0.09,0.1)\}$ for the ER graphs, $r$ values from $\{(0,0.02), (0.04,0.05), (0.08,0.12), (0.15,0.175)\}$ for GR graphs, and $k$ values from $\{(4,7), (12,18), (22,25), (27,30)\}$ for PA graphs. We refer to this collection as the {\em interval-drawn collection}.}

\blue{With the point-drawn collection, we aim to understand the relationships of the topological pseudometrics to the model parameter and to be able to determine if strong relationships between a given pseudometric and a topological metric may have originated from the latent parameters of the models. The interval point-drawn collection allows us to study how well a pseudometric predicts topological features. We chose model parameters that ensure that our datasets include graphs with genuinely different topologies and graphs with non-trivial $6$th Betti numbers.}

\subsection{Computing \blue{Topological Pseudometrics $d_\beta$ and $d_\Delta$}}
\label{sec:baselines}
For each test graph, we run our experiments by comparing well-established pseudometrics (TriadEuclid, TriadEMD, and PD) against the topological pseudometrics $d_\beta$ and $d_\Delta$ defined in Section~\ref{sec:prelim}.  

\para{Computing $d_\beta$ and $d_\Delta$ } 
To compute Betti numbers of directed graphs, we use the \textsf{Flagser}~\cite{LuetgehetmannGovcSmith2019} software available via a \textsf{Python} implementation \textsf{pyflagser}~\cite{TauzinLupoTunstall2020}. 
For homology computations, \textsf{Flagser} comes with an approximation option, which speeds up computation time while maintaining remarkable accuracy~\cite{LuetgehetmannGovcSmith2019}. 
In particular, it can be used to approximate the homology of a directed flag complex (associated with a test graph), by skipping columns that require more reduction steps than a user-chosen approximation parameter (denoted $\epsilon$). 

We use \textsf{Flagser} to compute the Betti numbers of the directed flag complex of each test graph up to dimension 6. \blue{We also compute approximate Betti numbers  with $\epsilon=10^0,10^1,10^2,10^3$, and $10^4$ to construct an approximate $d_{\beta}$ pseudometric.}  The corresponding errors of the logarithm of the Betti numbers (i.e., vectors of $max\{0, log(\beta_k)\}$) across our full data set is $36.1\%, 7.72\%, 1.53\%,  0.22\%$ and $0.04\%$, respectively. \blue{Since there is no theoretical error estimate for those approximate parameters, we calculated exact Betti numbers and then compared the approximated Betti numbers in order to calculate the errors above. We remark that similar calculations with $\epsilon = 10^5$ produced the exact Betti numbers.}

As explained in Section~\ref{sec:prelim}, we then compute the distance $d_\beta(G, G')$ between two directed graphs $G$ and $G'$ as the Euclidean distance between the vector of (possibly approximated) Betti numbers associated to the graphs $G$ and $G'$. 
The pseudometric $d_\Delta$ is defined similarly, with Betti numbers replaced by simplex counts.

\para{A parameter distance $d_p$ } 
In addition to $d_\beta$ and $d_\Delta$, we compute an additional \blue{pseudometric $d_p$ for comparison}, defined in terms of the parameters associated with our random graph models. 
For test graphs $G$ and $G'$ generated by a fixed random graph model (ER, GR, or PA), we define $d_p(G,G')$ to be the absolute value of the difference of the parameters used to generate them. 
For instance, if $G$ and $G'$ are generated with a ER random graph model with $\rho = 0.03$ and $0.06$ respectively, then $d_p(G, G') = |0.03-0.06|=0.03$. 

\para{Details of implementation } 
All calculations are performed using \textsf{Python 3.9}. The source code for our experiments can be found in our github repository \url{https://github.com/winsy17/graph_pseudo_top_view}.
Each ER/GR/PA graph is generated using the \textsf{Python} package \textsf{NetworkX}.

In the case of the point-drawn collection, we \blue{calculate the exact Betti numbers for all the parameters, except for $\rho \in \{0.15,0.2,0.25\}$, $r = 0.3$ and $k=70$ for which we use the approximation option of \textsf{pyflagser} with $\epsilon = 10^5$ because they are not computationally feasible to calculate exactly.} For the interval-drawn collection, we restrict the parameters to values that are computationally feasible and explore in \textsf{pyflagser}. 
We perform hierarchical clustering with complete linkage as our clustering technique. 
As a control, we also generate a random positive definite matrix with zeros on the diagonal that we include as a ``random'' pseudometric, denoted as ``random". 

\begin{figure}[!h]
    \includegraphics[width=0.98\linewidth]{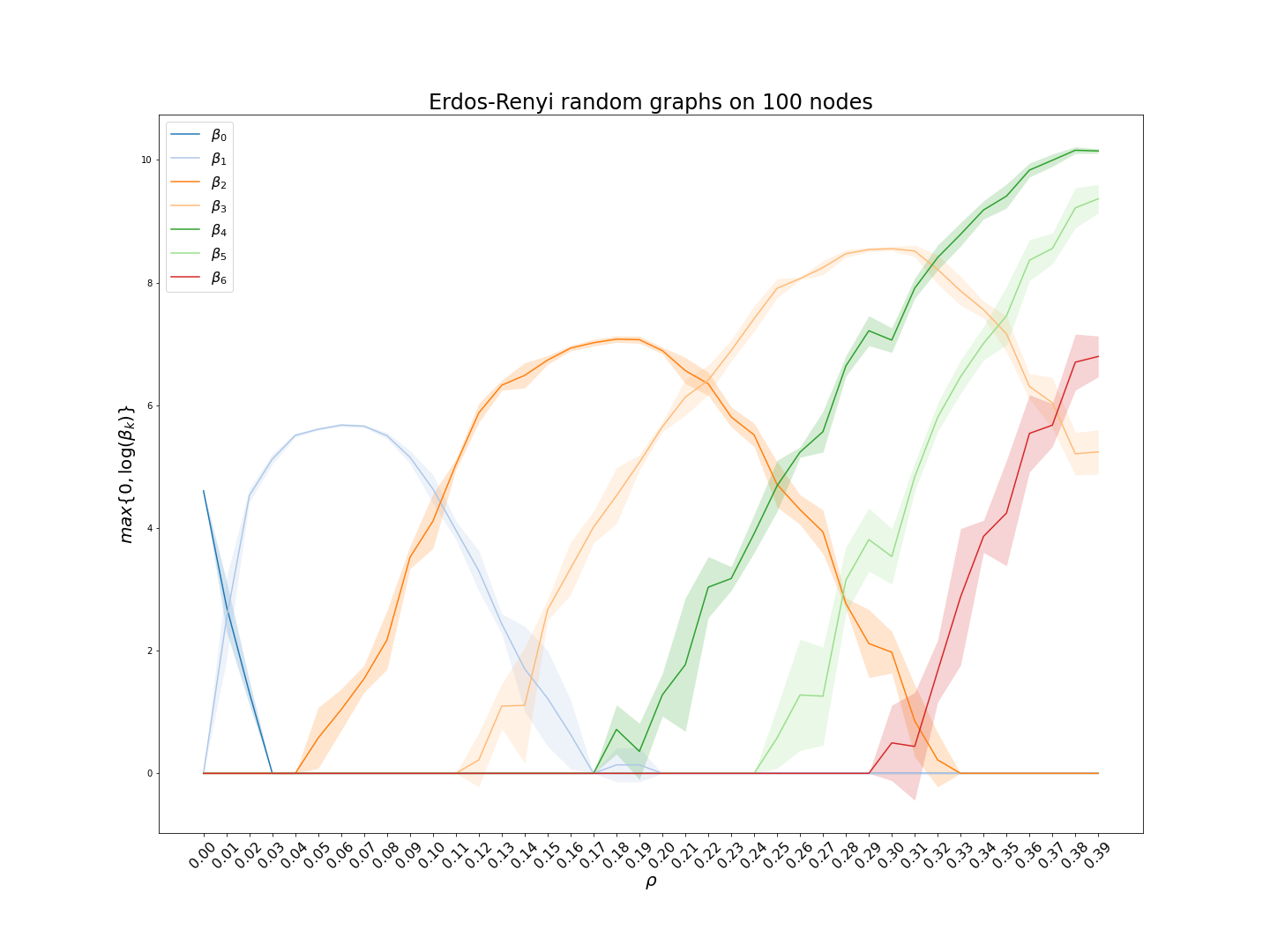}
    \vspace{-6mm}
    \caption{The logarithm of Betti numbers across a range of parameters for  Erdős–Rényi random graphs, each with $100$ vertices. The $x$-axis corresponds to the parameter $\rho$ ranging from $0$ to $0.39$, while the $y$-axis corresponds to the maximum of the logarithm of Betti numbers and 0. \blue{Similar results are observed for geometric random and preferential attachment random graphs.} }
    \label{fig:log_Betti_ER}
\end{figure}

\para{Distributions of Betti numbers }
As illustrated in~\autoref{fig:log_Betti_ER}, the Betti numbers seem to be strongly related to the parameters used to generate the test graphs.  
An analogous figure appears in \cite{KahleMeckesothers2013} for the undirected case. 
This is well aligned with the observations in \cite{Lasalle2019}, where Lasalle studied the behavior of Betti numbers for directed ER graphs.  

\section{Experimental Results}
\label{sec:results}

In this section, we first study the similarities between clusterings based on $d_\beta$ and $d_\Delta$ and those based on the three well-established pseudometrics by computing the Fowlkes-Mallows indices and distance correlations between them (Section~\ref{sec:FM-DC})\blue{, first for each of the three random graph models individually, then for all three combined. This analysis shows that when we consider them model by model, the pseudometrics we study are closely related to one another, with a few exceptions, which vary between models.   When all three models are considered simultaneously, the differences between the pseudometrics are substantially more stark.}

We next apply $k$-NN classification and report the classification accuracy for each of the three models of random graphs by model parameter \blue{and for the full collection of all of our random graphs by model type (Section~\ref{sec:accuracy}).  We achieve 100\% accuracy in almost all cases for the point-drawn collection.  In the interval-drawn case, MSE is very low for the ER and GR models and rather high for the PA model.}

We then explore the relationship of each of our pseudometrics to the various random graph parameters  by performing permutation tests with distance correlation \blue{between a fixed pseudometric and $d_\beta$ or $d_\Delta$, for datasets from the interval-drawn collections (Section~\ref{sec:permutation}).  When we consider all graph models together, the very small $p$-value enables us to reject the null hypothesis of independence of the pseudometrics.  The results are more nuanced when we instead study the models individually, one parameter at a time, where we can reject the null hypothesis of independence only rarely and then almost only for the relationship between $d_\beta$ and $d_\Delta$.}

Finally, we use $k$-NN regression on our pseudometrics \blue{in the interval-drawn case} to predict the topological feature vectors $b(G)$ or $c(G)$ and report the MSE (Section~\ref{sec:predict}).  \blue{Both $d_\beta$ and $d_\Delta$ perform well in predicting $b(G)$ and $c(G)$, better than any other pseudometric for the GR and PA models.  The case of the ER model is more complex, as (unsurprisingly) $d_\Delta$ and $d_\beta$ best predict $c(G)$ and $b(G)$, respectively, but the other pseudometrics perform quite well, too.}


\subsection{Fowlkes-Mallows Indices and Distance Correlation Between Pseudometrics}
\label{sec:FM-DC}

\begin{figure}[h]
    \centering
    \includegraphics[width=\textwidth]{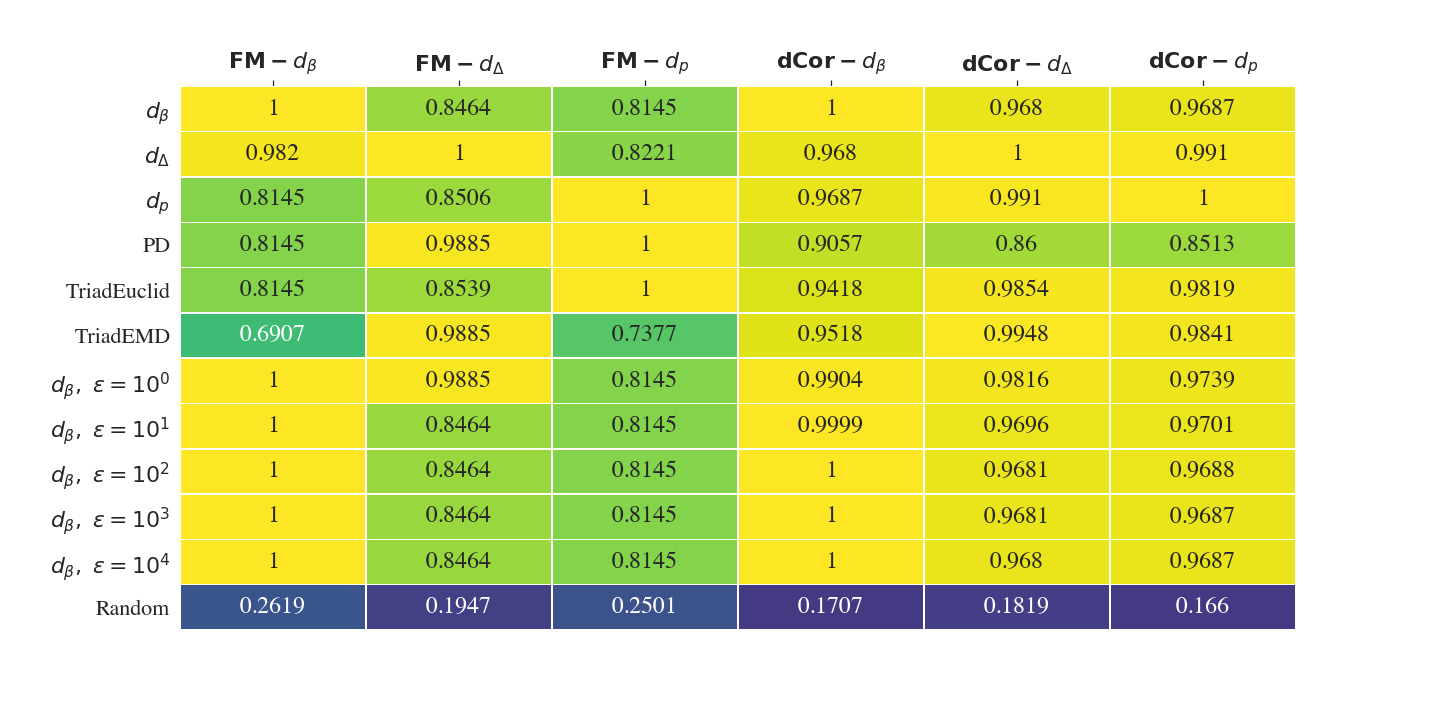}
    \vspace{-1cm}
    \caption{Interval-drawn collection: Erdős–Rényi (ER) random graphs}
\label{table:ER-FM-DC}
\end{figure}

\begin{figure}[h]
    \centering
    \includegraphics[width=\textwidth]{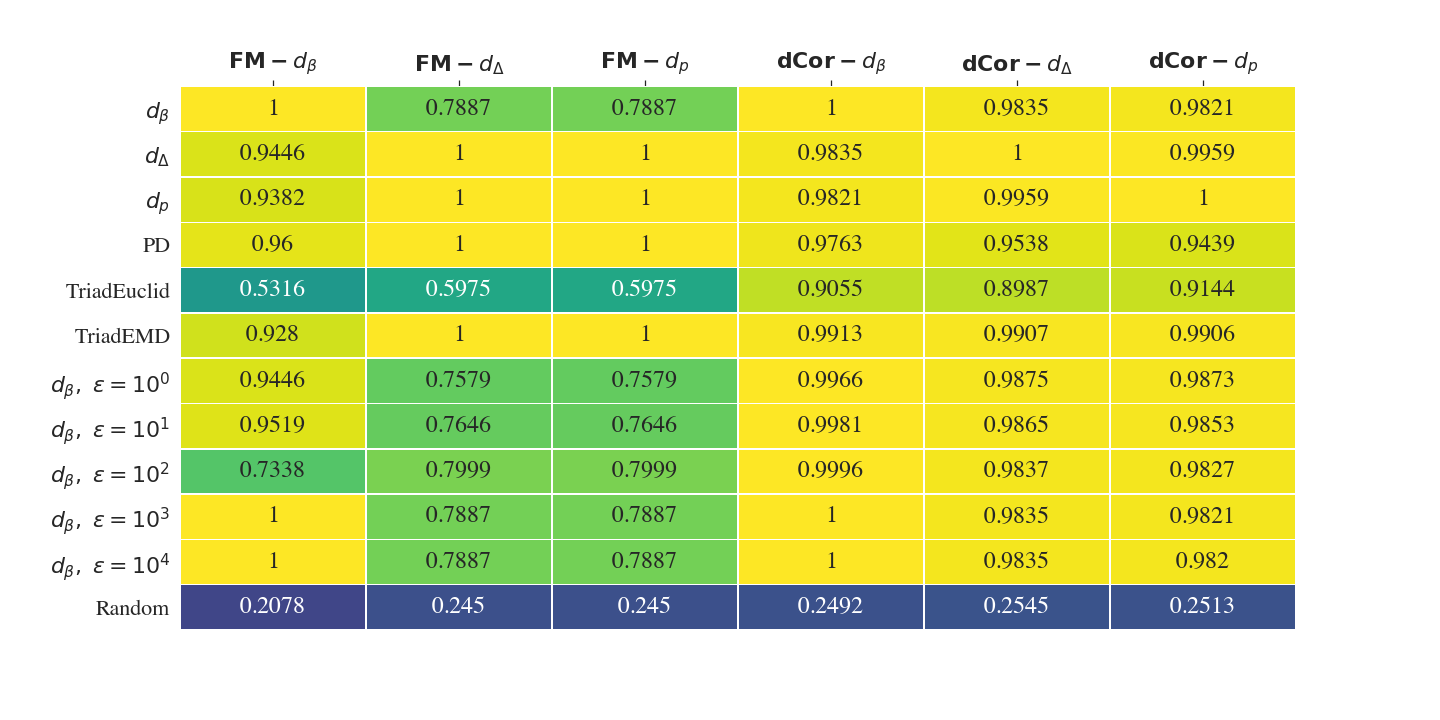}
    \vspace{-1cm}
\caption{Interval-drawn collection: geometric random (GR) graphs} 
\label{table:GR-FM-DC}
\end{figure}

\begin{figure}[h]
    \centering
    \includegraphics[width=\textwidth]{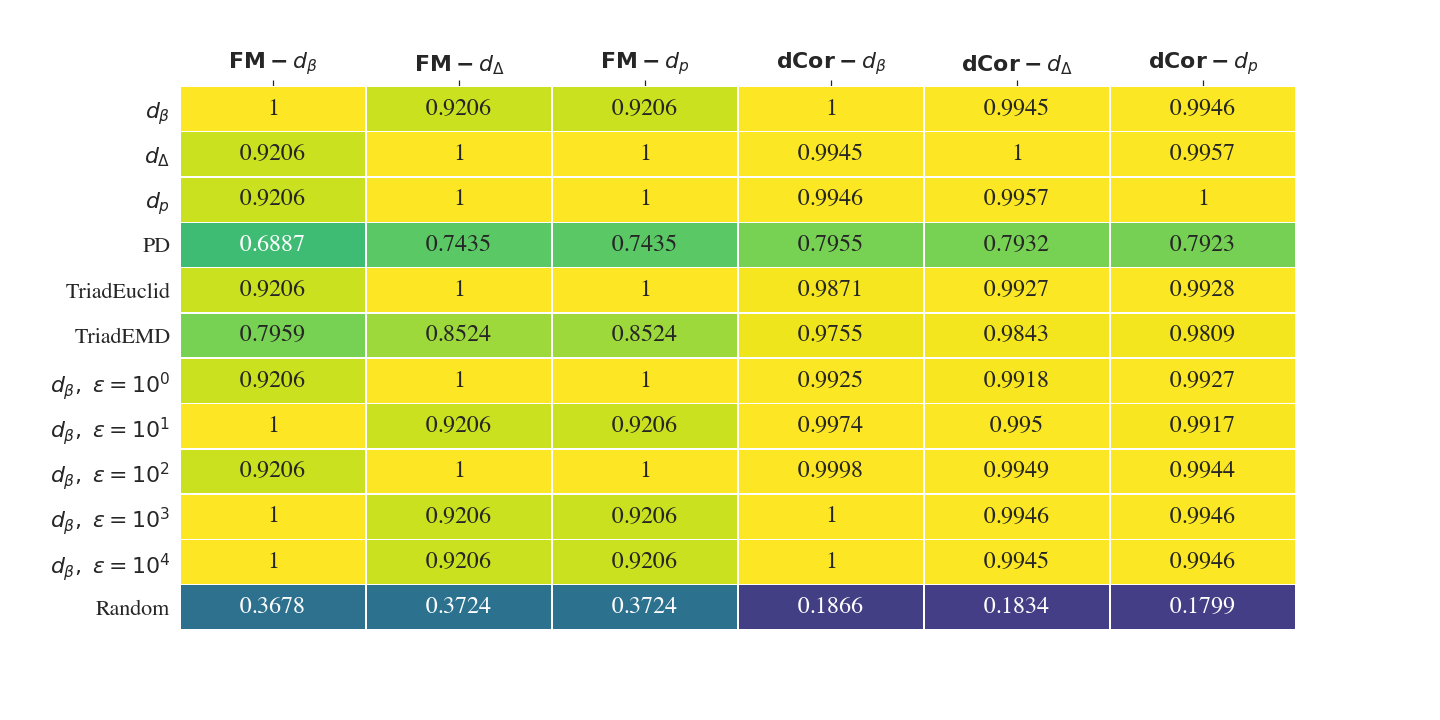}
    \vspace{-1cm}
\caption{Interval-drawn collection: preferential attachment (PA) random graphs} 
\label{table:PA-FM-DC}
\end{figure}

\begin{figure}[h]
    \centering
    \includegraphics[width=\textwidth]{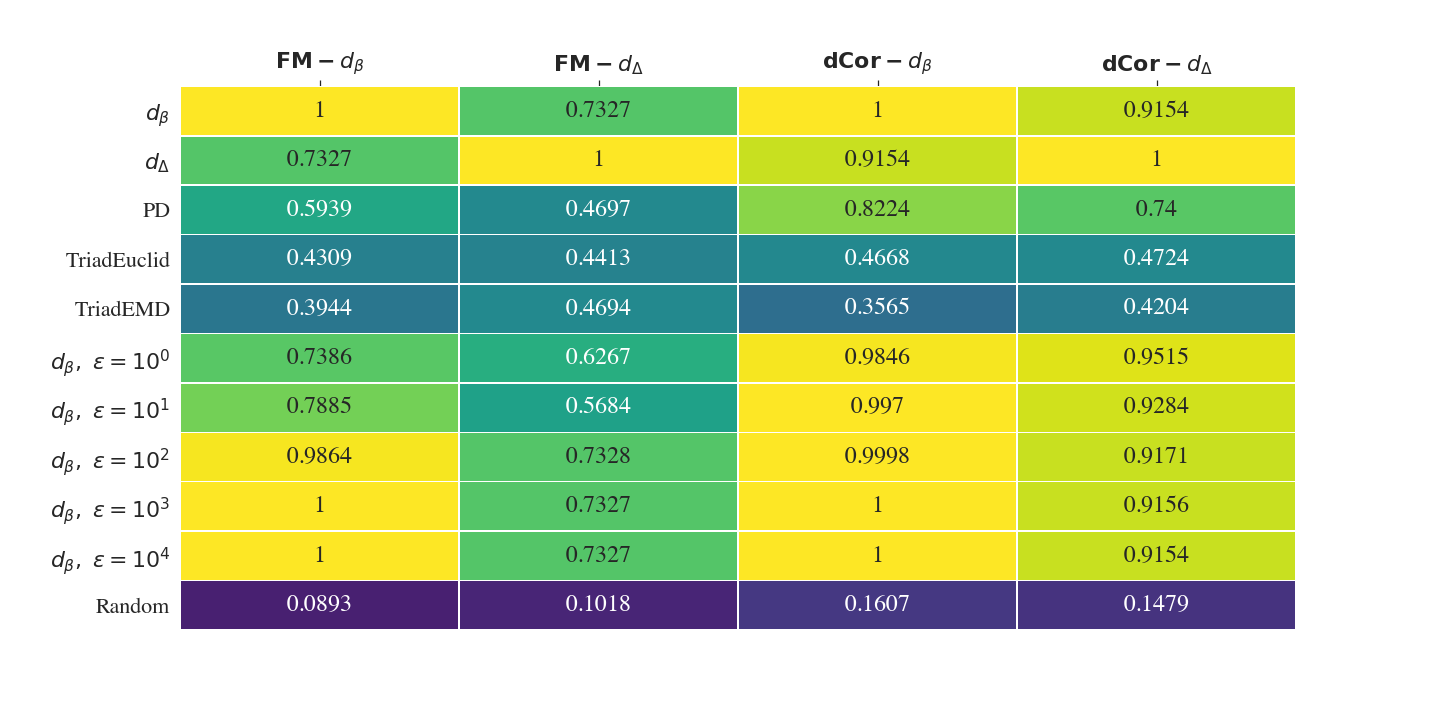}
    \vspace{-1cm}
\caption{Interval-drawn collection with all three random graph models} 
\label{table:all-FM-DC}
\end{figure}

We compute the Fowlkes-Mallows (FM) indices and distance correlations  of various pairs of  pseudometrics for the interval-drawn collection of random graphs; recall that this contains 100 graphs for each of the three models, ER, GR, and PA. 
We first treat each set of 100 graphs separately and report the resulting FM indices and distance correlations for ER, GR, and PA in~\autoref{table:ER-FM-DC}, \autoref{table:GR-FM-DC}, and \autoref{table:PA-FM-DC}, respectively. 
We then combine all 300 graphs and report the results in~\autoref{table:all-FM-DC}. 

Recall that a higher value for the FM index indicates greater  similarity between the clusterings induced by two pseudometrics. 
On the other hand, a higher value of distance correlation measures a higher level of dependence between the two pseudometrics. Both the FM index and distance correlation take values between 0 and 1.
 
Starting from~\autoref{table:ER-FM-DC}, we study various pseudometrics for  a collection of 100 ER random graphs. \blue{
The first column lists} the pseudometrics we consider: from top to bottom, the topological pseudometrics based on Betti numbers ($d_\beta$) and simplex counts ($d_\Delta$);  the  pseudometric $d_p$ based on model parameters; well-established (previously known) pseudometrics including PD (Portrait Divergence), TriadEuclid, and TriadEMD; $d_\beta$ using approximated Betti numbers with approximating parameters ($\epsilon = 10^0, 10^1, 10^2, 10^3, 10^4$); and a 	``random'' pseudometric generated with random positive definite matrix with zeros in the diagonal.  
The column labeled ``\textbf{FM-$d_\beta$}'' computes the FM index between $d_\beta$ and the pseudometric of each row (treated as the benchmark classification). 
For instance, the FM index between $d_{\beta}$ and $PD$ is $0.8145$. 
The number of clusters used for the FM index computation is chosen by a silhouette analysis of the column pseudometric. \blue{Notice that this way of choosing the number of clusters may result in different FM indices of the same two pseudometrics. For example, in \autoref{table:ER-FM-DC} the FM index of $d_\Delta$ with respect to $d_\beta$ differs to the FM index of $d_\beta$ with respect to $d_\Delta$.}
Similarly, the column ``\textbf{dCor-$d_\Delta$}'' computes the sample distance correlation between $d_\Delta$ and a pseudometric of each row. 
For instance, the sample distance correlation between $d_\Delta$ and PD is $0.86$. 

Considering ~\autoref{table:ER-FM-DC} closely, we notice that almost all the pseudometrics (except the random metric in the last row) are closely related, since the table contains high FM indices and high distance correlations. For simplicity, from now on, we exclude the last row (random) from our discussion. 

In particular, the \textbf{FM-$d_\beta$} column of ~\autoref{table:ER-FM-DC} indicates that the clusters based on $d_\beta$ are very similar to those based on $d_\Delta$ (treated as the benchmark classification), with an FM index of $\textbf{0.982}$.  
Moreover, $d_\beta$ is highly correlated with $d_\Delta$, with a distance correlation of $\textbf{0.968}$ in the \textbf{dCor-$d_\beta$} column. 
In addition, $d_\beta$ is also highly correlated with $d_p$ where $dCor(d_\beta, d_p) = \textbf{0.9687}$. 
However, as the TriadEMD entries in the \textbf{FM-$d_\beta$} and \textbf{dCor-$d_\beta$} columns show, high distance correlation between metrics ($d_\beta$ vs. TriadEMD, $\textbf{0.9518}$) does not imply similar clusterings:  the clusters given by $d_\beta$ are less similar to those given by TriadEMD than for any other pseudometric, with an FM index of $\textbf{0.6907}$. 

\blue{We remark that it is not too surprising that some of the numbers of ~\autoref{table:ER-FM-DC} associated to the Fawlkes-Mallows index coincide. For example, this occurs for \textbf{FM-$d_\beta$} of the pseudometrics $d_p$, PD and TriadEuclid. This phenomenon is due to identical clustering of the pseudometrics $d_p$, PD and TriadEuclid at the number of clusters chosen for $d_{\beta}$. This phenomenon can also be observed in subsequent figures/tables.} 

Careful inspection of~\autoref{table:GR-FM-DC} and~\autoref{table:PA-FM-DC} leads to a similar observation that the pseudometrics we study are highly related to one another (with a few exceptions -- with FM index less than $0.7$).  
In particular, $d_\beta$ and $d_\Delta$ consistently have high FM indices and high distance correlations across all three models. 
There are also some marked differences among the three models. 
For instance, the FM index between $d_{\beta}$ and PD for PA random graphs is quite low --  $\textbf{0.6887}$ in the \textbf{FM-$d_\beta$} column of \autoref{table:PA-FM-DC}  -- in comparison with the other two models. 

We combine all 300 random graphs  to study FM indices and distance correlations in~\autoref{table:all-FM-DC}. 
Both $d_\beta$ and $d_\Delta$ remain reasonably similar to and dependent on one another, with an FM index of $0.7327$ and a distance correlation of $0.9154$. 
We also see a sharp fall in the FM indices and distance correlations between the the topological ($d_\beta$, $d_\Delta$) and well-established pseudometrics (PD, TriadEuclid, TriadEMD). 
One potential explanation for these low FM indices and distance correlations  is that the relationship observed model by model is due to  a latent variable of the chosen model of random graph in addition to the parameter of that model.

\subsection{Classification Accuracy}
\label{sec:accuracy}

Recall that our point-drawn collection of random graphs consists of 10 ER random graphs generated for each of six parameters (six labels), 10 GR graphs generated for each of three parameters (three labels), and 10 PA graphs generated for each of three parameters (three labels). 
We apply $k$-NN classification and report the classification accuracy for each of the three  random graph models, as presented in \autoref{table:classify-three}, \autoref{table:classify-combined}, and \autoref{table:regression-three}.
 
As shown in~\autoref{table:classify-three}, using the PD, $d_\beta$, or $d_\Delta$  pseudometrics, we achieve 100\% accuracy for all three random graph models; the accuracy for TriadEuclid and TriadEMD is slightly lower.  

\begin{table}[H]
 \centering 
 \begin{tabular}{c| c c c}
\textbf{Pseudometrics} & \textbf{ER} & \textbf{GR} & \textbf{PA}  \\ [0.5ex]
\hline 
PD& 1 & 1 & 1\\ 
TriadEuclid & 0.9167 & 1 & 1\\ 
TriadEMD & 0.9167 & 1 & 1\\ 
$d_\beta$ & 1 & 1 & 1\\ 
$d_\Delta$ & 1 & 1 & 1\\ 
Random & 0.1833 & 0.2667 & 0.2\\ 
\end{tabular}\caption{Point-drawn collection: parameter classification rate within each set of random graphs}
\label{table:classify-three}
 \end{table}
 
If we treat all 120 random graphs in the point-drawn collection together and try to classify them into 3 classes (ER, GR and PA), we achieve 100\% classification accuracy for all pseudometrics (excluding the random metric), as shown in~\autoref{table:classify-combined}. 

\begin{table}[H]
 \centering 
 \begin{tabular}{c| c c c}
\textbf{Pseudometrics} & \textbf{Classification rate}   \\ [0.5ex]
\hline 
PD& 1\\ 
TriadEuclid & 1\\ 
TriadEMD & 1\\ 
$d_\beta$ & 1\\ 
$d_\Delta$ & 1\\ 
Random & 0.3917\\ 
\end{tabular}\caption{Point-drawn collection: model classification rate for all three sets of random graphs combined}
\label{table:classify-combined}
 \end{table}
 
For the interval-drawn collection, we apply $k$-NN regression and report the MSE in predicting the model parameters. 
As shown in~\autoref{table:regression-three}, for both ER and GR random graphs, the regression using both well-established and topological pseudometrics achieves very low MSE (excluding the random metric); while for the PA random graphs, none of these pseudometrics performs well.   

\begin{table}[H]
 \centering 
 \begin{tabular}{c | c c c}
\textbf{Pseudometric} & \textbf{ER} $(\times e-06)$ & \textbf{GR} $(\times e-06)$ & \textbf{PA}  \\ [0.5ex]
\hline 
PD& 0.9787 & 7.269 & 0.1488\\ 
TriadEuclid & 0.8367 & 447.6 & 0.1656\\ 
TriadEMD & 2.207 & 11.83 & 0.6636\\ 
$d_\beta$ & 2.956 & 6.665 & 0.4048\\ 
$d_\Delta$ & 1.653 & 9.03 & 0.4084\\ 
Random & 1373 & 4047 & 98.3\\ 
\end{tabular}
\caption{Interval-drawn collection: MSE with $k$-NN regression in predicting model  parameters}
\label{table:regression-three}
 \end{table}
 

\subsection{Permutation Tests}
\label{sec:permutation}

We are interested in the statistical significance of relationships between well-established pseudometrics and $d_\beta$ (respectively, $d_\Delta$). Permutation tests provide a method to test statistical significance of high distance correlation or Fowlkes-Mallows (FM) index as demonstrating dependence between two pseudometrics. Our null hypothesis is that the two pseudometrics being compared are independent. 

\blue{
For small $p$-values we can reject the null hypothesis and deduce dependency. 
For insufficiently small $p$-value we cannot reject the null hypothesis of independence. This does not imply the two pseudometrics are independent, but rather that we cannot rule out the independence.}

\blue{We perform permutation tests across all four datasets from the interval-drawn collection, based on either FM index or distance correlation as the dependency measures. 
One of our null-hypotheses is thus the independence of pseudometrics from our topological metrics ($d_\beta$, $d_\Delta$ or $d_p$) as measured by the FM index. 
We also consider the analogous null-hypothesis of independence with respect to distance correlation.} 


\blue{With exception of the random pseudometric, all of the calculated $p$-values equal $0.0006667$. Hence we can reject both of our null-hypotheses with $p$-value equal to $0.0006667$. }

\blue{Except for the random control, the FM indices (respectively, distance correlations) calculated in Section~\ref{sec:FM-DC} are so high that for every single random permutation the result is lower than for the original (unpermuted) FM indices (respectively, distance correlations). With $1500$ permutations used in these tests, we compute a $p$-value of $0.000667$. In each case we have performed at most $66$ tests. Since $0.000667<0.05/66$ we reject the null-hypothesis (except the random control) with a family-wise error of $0.05$.} 

\blue{If we were to consider all of the permutation tests combined, then the number of tests ($264$) would require a prohibitive number of permutations to make the Bonferroni correction useful. We would require $p$-values of at most $0.000189$ to be able to establish any test as significant, which would require a minimum of $5280$ permutations for each test.}

\blue{We use instead the Benjamini-Hochberg procedure to combine all the tests. Since $$0.000667<\frac{242}{262*C(262)}{0.005},$$  we can reject all the null hypotheses for all the pseudometrics except the random controls and still have a false discovery rate bounded above by $0.005$. 
Here we have used that $C(264)<\ln(264)+1<6.57$. }

The results in Section~\ref{sec:FM-DC} indicated that all of the pseudometrics are related in terms of what they detect in our random graph models. However, it is unclear to what extent this relationship might also rely on or be amplified by the differences between the models and model parameters chosen. It is therefore necessary to test conditional independence while taking these possible effects into account.

One way to control for these latent variables is to use the samples with fixed choices of generating model and model parameter, that is, analysing slices of the point-drawn collection taken based on parameter for each model.
%
%
%
\autoref{table:p-ER-Betti} through \autoref{table:p-PA-simplex} contain the $p$-values of the permutation tests based on distance correlations between a pseudometric (in the first 1st column) and $d_\beta$ or $d_\Delta$. 

\blue{
We reject the null hypothesis test for independence when the $p$ value is $0.0004998$, which corresponds to the permutation tests where none of the random permutations had a smaller distance correlation than that of the original labels. These $p$-values are in bold.}

\blue{To control for multiple hypothesis testing, we again consider the family-wise error over each table separately and then the false discovery rate over the tables combined. Tables 5 to 10 each have at most $30$ entries, so we can say that all the results with values $p<0.0005$ are significant with a family-wise error bound of $0.015$. }

\blue{We can again use the the Benjamini-Hochberg procedure to bound the false discovery rate on the $p$-values written in bold in Tables 5 to 10 combined. Since $$0.0004998<\frac{7}{120*C(120)}{0.03},$$ (using $C(120)<\ln(120)+1<3.5$) we can report all of the $p$-values less than $0.0005$ as significant, with an expected  false discovery rate bounded above by $0.03$. This means the expected number of false positives is less that $0.21$. }

\blue{Most of the rejected null hypotheses concern the independence of distances based on Betti numbers and on simplex counts. Furthermore these are occurring with the higher parameter values within the models. This is not surprising, as it ties in well with the literature on limit theorems for Betti numbers (see the survey paper \cite{Bobrowski} and its references).}

\blue{What is more surprising is that for the geometric random graph with parameter value $0.3$, we can reject independence between the simplex count metric and the portrait divergence metric. This suggests that for geometric graphs, the differences in simplex counts is in fact more closely related to the portrait divergence than it is to differences in the Betti numbers. }

\blue{Although the various pseudometrics are not independent when we consider sets of directed graphs with multiple model parameters, we generally cannot reject independence once we restrict to specific models and parameter. This indicates that the relationships between these pseudometrics are probably driven by the latent variable of the models and their parameters. They may be capturing different features within these graphs (that are all affected in different ways by the model and parameter) and in doing so provide different perspectives on the graph structure.
}

\begin{table}[H]
 \centering 
 \begin{tabular}{c | c c c c c c}
Pseudometrics  & 0.03 & 0.06 & 0.1 & 0.15 & 0.2 & 0.25\\ [0.5ex]
\hline 
PD& 0.5442 & 0.7251 & 0.1529 & 0.3673 & 0.8041 & 0.01349\\ 
TriadEuclid & 0.1339 & 0.4803 & 0.8281 & 0.6887 & 0.2304 & 0.3318\\ 
TriadEMD & 0.6632 & 0.1964 & 0.7206 & 0.9245 & 0.2909 & 0.2539\\ 
$d_\Delta$ & 0.6387 & 0.7916 & 0.5922 & 0.5587 & 0.1399 & \textbf{0.0004998}\\ 
Random & 0.3603 & 0.5882 & 0.3693 & 0.2614 & 0.8431 & 0.3003\\ 
\end{tabular}
\caption{ER random graphs in point-drawn collection: $p$-values for distance correlation w.r.t. $d_\beta$}
\label{table:p-ER-Betti}
\end{table}

\begin{table}[H]
 \centering 
 \begin{tabular}{c | c c c c c c}
Pseudometrics  & 0.03 & 0.06 & 0.1 & 0.15 & 0.2 & 0.25\\ [0.5ex]
\hline 
PD& 0.96 & 0.5187 & 0.3018 & 0.3638 & 0.5927 & 0.01299\\ 
TriadEuclid & 0.7031 & 0.9535 & 0.6592 & 0.3973 & 0.8661 & 0.3188\\ 
TriadEMD & 0.3138 & 0.4123 & 0.9355 & 0.4478 & 0.988 & 0.1939\\ 
$d_\beta$ & 0.6622 & 0.8061 & 0.6047 & 0.5592 & 0.1489 & \textbf{0.0004998}\\ 
Random & 0.7056 & 0.6022 & 0.1084 & 0.5872 & 0.9455 & 0.2454\\ 
\end{tabular}
\caption{ER random graphs in point-drawn collection: $p$-values for distance correlations w.r.t. $d_\Delta$}
\label{table:p-ER-simplex}
 \end{table}

\begin{table}[H]
 \centering 
 \begin{tabular}{c | c c c }
Pseudometrics  & 0.1 & 0.175 & 0.3\\ [0.5ex]
\hline 
PD& 0.03348 & 0.2034 & 0.02299\\ 
TriadEuclid & 0.5642 & 0.2699 & 0.2519\\ 
TriadEMD & 0.1034 & 0.02049 & 0.906\\ 
$d_\Delta$ & 0.2769 & 0.2199 & 0.05597\\ 
Random & 0.01449 & 0.1099 & 0.6187\\ 
\end{tabular}
\caption{GR random graphs in point-drawn collection: $p$-values for distance correlation w.r.t. $d_\beta$}
\label{table:p-GR-Betti}
\end{table}

\begin{table}[H]
 \centering 
 \begin{tabular}{c | c c c }
Pseudometrics  & 0.1 & 0.175 & 0.3\\ [0.5ex]
\hline 
PD& 0.3903 & 0.02199 & \textbf{0.0004998}\\ 
TriadEuclid & 0.5312 & 0.7641 & 0.8221\\ 
TriadEMD & 0.7536 & 0.92 & 0.7866\\ 
$d_\beta$ & 0.2909 & 0.2214 & 0.06197\\ 
Random & 0.5557 & 0.3488 & 0.3158\\ 
\end{tabular}
\caption{GR random graphs in point-drawn collection: $p$-values for distance correlation w.r.t. $d_\Delta$}
\label{table:p-GR-simplex}
 \end{table}

\begin{table}[H]
 \centering 
 \begin{tabular}{c | c c c }
Pseudometrics  & 20 & 40 & 70\\ [0.5ex]
\hline 
PD& 0.4858 & 0.07796 & 0.4513\\ 
TriadEuclid & 0.8981 & 0.7466 & 0.4943\\ 
TriadEMD & 0.1904 & 0.1074 & 0.5762\\ 
$d_\Delta$ & 0.3058 & \textbf{0.0004998} & \textbf{0.0004998}\\ 
Random & 0.1979 & 0.7221 & 0.08796\\ 
\end{tabular}
\caption{PA random graphs in point-drawn collection: $p$-values for distance correlation w.r.t. $d_\beta$}
\label{table:p-PA-Betti}
 \end{table}

\begin{table}[H]
 \centering 
 \begin{tabular}{c | c c c}
Pseudometrics & 20 & 40 & 70\\ [0.5ex]
\hline 
PD& 0.947 & 0.04998 & 0.4253\\ 
TriadEuclid & 0.2029 & 0.6912 & 0.6617\\ 
TriadEMD & 0.9185 & 0.06047 & 0.6342\\ 
$d_\beta$ & 0.3013 & \textbf{0.0004998} & \textbf{0.0004998}\\ 
Random & 0.6522 & 0.3288 & 0.09545\\ 
\end{tabular}
\caption{PA random graphs in point-drawn collection: $p$-values for distance correlation w.r.t. $d_\Delta$}
\label{table:p-PA-simplex}
\end{table}

\subsection{Comparison of Clusterings and Classification Power}
\label{sec:predict}

We now turn to testing how well each pseudometric performs in classification/regression tasks with respect to our graph models and parameters. 
For this, we focus on the interval-drawn collection, applying $k$-NN regression with respect to our pseudometrics to predict the topological feature vectors $b(G)$ and $c(G)$. 
Note that we have chosen the model parameters deliberately so as to potentially confuse any Betti number classifier, so that we are being very conservative about the potential power of Betti numbers for classification in a mixed model scenario.
In the tables below, we use one of the pseudometrics of the first column to predict our topological feature vector and report the MSE results. 

For individual random graph models, it is not surprising that $d_\beta$ and $d_\Delta$ predict both simplex counts and Betti numbers fairly well. 
In the case of mixed model data in  (\autoref{table:predict-Betti-all} and \autoref{table:predict-simplex-all}), it is clear that PD outperforms the other well-established pseudometrics, as measured by MSE. \blue{In this same setting, the very large MSEs from the predictions using TriadEuclid and TriadEMD make them unsuitable methods for these classification tasks.}


\begin{table}[H]
 \centering 
 \begin{tabular}{c | c c c}
Pseudometric & ER & GR & PA  \\ [0.5ex]
\hline 
PD& 0.4179 & 0.6574 & 0.4238\\ 
TriadEuclid & 0.4969 & 3.569 & 0.4338\\ 
TriadEMD & 0.4504 & 0.8846 & 0.5576\\ 
$d_\beta$ & 0.1263 & 0.2526 & \textbf{0.08243}\\ 
$d_\Delta$ & 0.5233 & 0.588 & \textbf{0.3343}\\ 
$d_p$ & 0.4736 & 0.6567 & 0.421\\ 
Random & 46.15 & 39.86 & 25.92\\ 
\end{tabular}
\caption{Interval-drawn collection: MSE in predicting Betti numbers}
\label{table:predict-Betti}
 \end{table}

\begin{table}[H]
 \centering 
 \begin{tabular}{c | c c c}
Pseudometric & ER & GR & PA  \\ [0.5ex]
\hline 
PD& 0.3368 & 0.8528 & 0.4618\\ 
TriadEuclid & 0.3044 & 14.6 & 0.4803\\ 
TriadEMD & 0.4215 & 0.9439 & 0.4793\\ 
$d_\beta$ & 0.5266 & 0.7196 & \textbf{0.343}\\ 
$d_\Delta$ & 0.124 & 0.2695 & \textbf{0.08164}\\ 
$d_p$ & 0.2638 & 0.8092 & 0.3726\\ 
Random & 51.31 & 130.6 & 33.12\\ 
\end{tabular}
\caption{Interval-drawn collection: MSE in predicting simplex counts}
\label{table:predict-simplex}
 \end{table}

\begin{table}[H]
 \centering 
 \begin{tabular}{c | c}
Pseudometric & MSE   \\ [0.5ex]
\hline 
PD& \textbf{0.7495}\\ 
TriadEuclid & 77.17\\ 
TriadEMD & 77.45\\ 
$d_\beta$ & 0.1378\\ 
$d_\Delta$ & 0.4866\\ 
Random & 45.98\\ 
\end{tabular}
\caption{Interval-drawn collection: MSE in predicting Betti numbers, combining  all models}
\label{table:predict-Betti-all}

 \end{table}

\begin{table}[H]
 \centering 
 \begin{tabular}{c | c}
Pseudometric & MSE   \\ [0.5ex]
\hline 
PD& \textbf{2.53}\\ 
TriadEuclid & 91.72\\ 
TriadEMD & 92.47\\ 
$d_\beta$ & 1.183\\ 
$d_\Delta$ & 0.1407\\ 
Random & 79.16\\ 
\end{tabular}
\caption{Interval-drawn collection: MSE in predicting simplex count, combining  all models}
\label{table:predict-simplex-all}
 \end{table}

\section{Conclusion}
\label{sec:conclusion}

\blue{In recent years several useful pseudometrics have been defined on the set $\G$ of finite directed graphs (digraphs), as tools for quantifying similarities and differences between digraphs. In this preliminary study, inspired by recent successful applications of topological methods to digraphs, we introduced two ``topological'' pseudometrics on $\G$, $d_\beta$ and $d_\Delta$, derived from Betti numbers and simplex counts of the directed flag complex of a digraph, respectively. We compared these new pseudometrics with three well-established pseudometrics (PD, TriadEuclid, TriadEMD), to determine to what extent the latter might be sensitive to topological structure.} 

\blue{The relationship between the topological pseudometrics $d_\beta$ and $d_\Delta$ and the previously defined pseudometrics  proved to be highly model-dependent.} For example, we generally observed high values of Fowlkes-Mallows indices and distance correlations between all pairs of pseudometrics for a fixed random model (\autoref{table:ER-FM-DC}, \autoref{table:GR-FM-DC}, \autoref{table:PA-FM-DC}),  but significantly lower values when mixing models (\autoref{table:all-FM-DC}).

The relationship between $d_\beta$ and $d_\Delta$ themselves is worth exploring, as calculating simplex counts is considerably computationally cheaper than calculating Betti numbers. When restricting to a specific random model, we obtained Fowlkes-Mallows indices between $d_\Delta$ and $d_\beta$ greater than $\textbf{0.9206}$ (column 1 of \autoref{table:ER-FM-DC}, \autoref{table:GR-FM-DC}, \autoref{table:PA-FM-DC}). The scores were noticeably lower, however, when we mixed models (\autoref{table:all-FM-DC}), \blue{indicating that the relationship between $d_\beta$ and $d_\Delta$ is also model-dependent.} 

We also compared the performance of $d_\beta$ and $d_\Delta$ at two classification tasks for random graphs: estimating parameters within a random model and classifying graphs by random model (Section~\ref{sec:accuracy}). We observed that the performance of $d_\beta$ and $d_\Delta$ is comparable to that of PD, TriadEuclid, and TriadEMD in both of these tasks, achieving almost 100\% classification accuracy (see \autoref{table:classify-three}, \autoref{table:classify-combined}) or very low MSE (see \autoref{table:regression-three}), \blue{except in the case of the PA model. It would be interesting to determine the source of the relatively large MSE in this last case.}

We explored, moreover, whether the model parameter explains the observed relationships between our topological metrics and PD, TriadEMD and TriadEuclid, using permutation tests with distance correlation and a null-hypothesis of independence  (Section~\ref{sec:permutation}). In only very few cases, \blue{almost always concerning the pair ($d_\beta,d_\Delta)$}, were we able to reject the null-hypothesis, \blue{so that the independence of most pairs of pseudometrics considered remains an open question for the graph models we studied.} 

Finally, we tested how well each pseudometric predicts Betti numbers and simplex counts (Section~\ref{sec:predict}). With exception of TriadEuclid on GR graphs (\autoref{table:predict-simplex}), all of the pseudometrics performed very well with a very low MSE when predicting either Betti or simplex counts for a fixed random model (\autoref{table:predict-Betti}, \autoref{table:predict-simplex}). When combining random models, applying PD resulted in a very low MSE when predicting either Betti or simplex counts, while using either TriadEuclid or TriadEMD led to significant prediction error for both Betti and simplex counts.   \blue{On the other hand, both $d_\beta$ and $d_\Delta$ performed very well in the case of combined models. It follows in particular that when models underlying a collection of random graphs are either known and mixed or unknown, one should use either PD or one of $d_\beta$ and $d_\Delta$ to predict Betti numbers or simplex counts.  Applying $d_\Delta$ is slightly less accurate in its prediction of Betti numbers than $d_\beta$, which should be weighed against its ease of computation.}

Although we considered only the topology of the directed flag complex in this study, our methods can be applied in a straightforward manner using other topological features, such as those obtained from the flag complex of an undirected graph or from the path homology of a directed graph~\cite{ChowdhuryMemoli2018, GrigoryanLinMuranov2020}.   

One limitation of this study is that all of the graphs analyzed have 500 vertices.  It would be interesting to determine whether the relationships observed here hold for graphs of other sizes.  In particular,  understanding the asymptotic behavior of these relationships, as the number of vertices goes to infinity, should provide important insights. 
The Betti numbers of various types of random simplicial complexes~\cite{KahleMeckesothers2013} and random ER graphs~\cite{Lasalle2019} have been studied previously. A similar limit theorem for Betti numbers of directed flag complexes seems possible, though nontrivial, to formulate and prove, given the computations displayed in ~\autoref{fig:log_Betti_ER}.


\begin{acknowledgement}
This work began at the Women in Computational Topology Workshop in 2019. The authors wish to thank the Australian Mathematical Sciences Institute, the National Science Foundation and the Association of Women in Mathematics for the accommodation and travel support. 
ALGP is supported by the EPSRC grant ``New Approaches to Data Science: Application Driven Topological Data Analysis'' EP/R018472/1. 
NY is supported by the Alan Turing Institute under the EPSRC grant EP/N510129/1. 
KT is supported by an ARC DECRA fellowship.
BW is supported in part by DOE DE-SC0021015 and NSF DBI-1661375.  
The authors would like to thank Erika Rold\'an for insightful contributions, Gesine Reinert for sharing the source code for TriadEMD and TriadEuclid, the Oxford Mathematical Institute for providing access to computational resources, and the anonymous referees whose comments have helped to clarify this paper.
\end{acknowledgement}


\bibliographystyle{abbrv}
\bibliography{directed-graphs-refs}

\end{document}